\documentclass[leqno,12pt]{amsart}
\usepackage{amssymb} 
\theoremstyle{plain} %text of this environment is typesetted in italics
\newtheorem{theorem}{\indent\sc Theorem}[section]
\newtheorem{lemma}[theorem]{\indent\sc Lemma}

\newtheorem{proposition}[theorem]{\indent\sc Proposition}

\theoremstyle{definition} 
\newtheorem{definition}[theorem]{\indent\sc Definition}
\newtheorem{remark}[theorem]{\indent\sc Remark}

\begin{document}

\title[$K$-theory of complete toric varieties]{$K$-theory of smooth complete toric 
varieties and related spaces}

\author[P. Sankaran]{Parameswaran Sankaran} 
%%%%%%%%%%%%%%% footnote %%%%%%%%%%%%%%%%


\address{%Parameswaran Sankaran
Institute of Mathematical Sciences \endgraf
CIT Campus, Taramani\endgraf
Chennai 600113 \endgraf
India}
\email{sankaran@imsc.res.in}

\maketitle

\def\theequation {\arabic{section}.\arabic{equation}}
\renewcommand{\thefootnote}{}
\newcommand{\codim}{\mbox{{\rm codim}$\,$}}
\newcommand{\stab}{\mbox{{\rm stab}$\,$}}

\newcommand{\lr}{\mbox{$\longrightarrow$}}
\newcommand{\bro}{{\Bbb R}_{\geq 0}}
\newcommand{\ct}{{\cal T}}
\newcommand{\cm}{{\cal M}}
\def\cO{\mathcal{O}}
\newcommand{\cd}{{\cal D}}
%\newsymbol \rightarrowtail 131A
\newcommand{\blr}{\Big \longrightarrow}
\newcommand{\da}{\Big \downarrow}
\newcommand{\ua}{\Big \uparrow}
\newcommand{\hra}{\mbox{\LARGE{$\hookrightarrow$}}}
\newcommand{\rt}{\mbox{\Large{$\rightarrowtail$}}}

\newcommand{\dua}{\begin{array}[t]{c}
\Big\uparrow \\ [-4mm]
\scriptscriptstyle \wedge \end{array}}
%%%%%%%%%%%%%%%  My macros %%%%%%%%

\def\ci{\mathcal{I}}
\newcommand{\cp}{{\cal P}}

\newcommand{\cch}{{\rm ch}}
\newcommand{\cg}{{\cal G}}
\newcommand{\bs}{{\Bbb S}}
\newcommand{\bz}{{\Bbb Z}}
\newcommand{\bc}{{\Bbb C}}
\newcommand{\ev}{{\rm ev}}
\newcommand{\bq}{{\Bbb  Q}}
\newcommand{\bzo}{{\Bbb Z}_{\geq 0}}
\newcommand{\bt}{{\Bbb T}}
\newcommand{\br}{{\Bbb R}}
\newcommand{\bp}{{\Bbb P}}
\newcommand{\wt}{\widetilde}

\newcommand{\ffi}{\frak{I}}
\newcommand{\fj}{\frak{J}}
\newcommand{\im}{{\rm Im}\,}
\newcommand{\gr}{{\rm gr}}
\newcommand{\inn}{{\rm in}}
\newcommand{\Hom}{{\rm Hom}}
\newcommand{\spin}{{\rm Spin}\,}
\newcommand{\ds}{\displaystyle}
\newcommand{\tor}{{\rm Tor}\,}
\def\ck{\mathcal{K}}
\def\ns{\mathop{\lr}}
\def\nssup{\mathop{\lr\,sup}}
\def\nsinf{\mathop{\lr\,inf}}
\newcommand{\eps}{\epsilon}

\newcommand{\bn}{{\Bbb N}}

\begin{abstract}
The $K$-rings of non-singular complex projective varieties as well as 
quasi-toric manifolds  
were described in terms of generators and relations in an earlier work of 
the author with V. Uma. 
In this paper we obtain a similar description for complete 
non-singular toric varieties. Indeed, our approach enables us to 
obtain such a description for the more general class of torus 
manifolds with locally standard torus action and orbit space
 a homology polytope. 
\end{abstract}

\section{Introduction}
This paper consists of two parts, the first of which gives a description 
of $K$-ring of a non-singular complete toric variety in terms of 
generators and relations. In the second part, which subsumes the 
first, we obtain the same result for the class of torus manifolds 
with locally standard action whose orbit space is a homology 
polytope.  Although the proofs in both cases involve the 
same steps, they are  easier to establish or well known in the 
context of toric varieties and helps to focus ideas. Besides, the 
language used in the two proofs are different. For these 
reasons, we have separated out the algebraic geometric case 
from the topological one, taking care to avoid unnecessary 
repetitions.  

Let $N\cong \bz^n$ and let $N_\br:=N\otimes_\bz\br$. We shall denote by $M$ 
the dual lattice $\Hom_\bz(N,\bz)$ and by $\bt$ the $n$-dimensional 
complex algebraic torus whose coordinate ring is the group algebra $\bc[M]$. 
The lattice $M$ is identified with the group of characters of $\bt$ 
and $N$ with the group of $1$-parameter subgroups of $\bt$. Let $\Delta$ be 
a fan in $N$ and let $X(\Delta)$ the corresponding complex $T$-toric variety.
When $\Delta$ is complete the variety $X(\Delta)$ is complete. When 
$\Delta$ is regular, (i.e the cones of $\Delta$ are generated by 
part of a $\bz$-basis for $N$) the variety $X(\Delta)$ is smooth. 
We refer the reader to \cite{fulton} for an introduction to 
toric varieties. 

Let $\Delta$ be a complete regular fan in $N$ and let $X:=X(\Delta)$ 
be the corresponding non-singular complete variety. 
Recall that one has an inclusion-reversing correspondence between 
cones of $\Delta$ and $\bt$-orbit closures of $X$. In fact the 
orbit closure $V(\sigma)$ of a cone $\sigma\in \Delta$ equals the 
union of $\bt$-orbits of a certain (closed) points $p_\gamma \in X$ as 
$\gamma$ varies over those cones of $\Delta$ which contain $\sigma$.
In particular, if $\sigma, \sigma'\in\Delta$ are not faces of a 
cone $\gamma\in\Delta$, then 
$V(\sigma)\cap V(\sigma')=\emptyset$.  The dimension of $V(\sigma)$ 
equals the co-dimension of $\sigma$ in $N_\br$.  The orbit closures 
are themselves toric varieties (under the action of 
suitable quotients of $\bt$) and are smooth as we assumed $X$ to be smooth.

We shall denote by $\Delta(k)$ the set of all $k$-dimensional cones of 
$\Delta$.

We shall now recall the description of the Chow ring $A^*(X)$ of $X$. 
The ring $A^*(X)$ is generated by the classes of divisors $[V(\rho)]$ 
as $\rho$ varies over the edges of $\Delta$. Also one has the following 
relations among these classes:\\

(i) if $\rho_1,\cdots ,\rho_k\in \Delta(1)$ 
do not span a cone in $\Delta$, then 
$$[V(\rho_1)]\cdots [V(\rho_k)]=0 \eqno{(1)}$$ 

(ii) For any $u\in M$ and $v_\rho\in N$ is the primitive vector on $\rho\in \Delta$ 
one has the relation 
$$\sum_{\rho\in \Delta(1)}\langle u,v_\rho\rangle[V(\rho)]=0.\eqno{(2)}$$ 

Danilov's theorem \cite[Theorem 10.8]{danilov} asserts that these are no further generating relations and that $A^*(X)$ has no torsion.
It turns out that the cycle class map  
$A^*(X)\lr H^*(X)$ is an isomorphism  
that doubles the gradation. So Danilov's result also yields a description 
of the singular cohomology (with $\bz$-coefficients).  In the case of 
non-singular projective toric varieties these results are due to 
Jurkiewicz \cite{jur}.

In this paper we give a description of the `topological' $K$-ring of $X(\Delta)$, 
denoted $K(X(\Delta))$ in terms of generators and relations. Denoting the 
Grothendieck ring 
of algebraic vector bundle  by $\ck(X(\Delta))$, it turns out 
that the forgetful homomorphism \newline
$\ck(X(\Delta)) \lr K(X(\Delta))$ is an 
isomorphism and hence we obtain a similar description of $\ck(X(\Delta))$ as well.
We shall also establish a similar description of the $K$-ring of a torus 
manifolds with locally standard action whose orbit space is a homology polytope.
These results were established for non-singular projective 
toric varieties in \cite{su1} and for  
quasi-toric manifolds in \cite{su2}. Recently V. Uma \cite{uma} has extended the 
results of \cite{su2} to the case torus manifolds under an additional 
shellability hypothesis on the orbit space which is assumed to be locally 
standard with quotient a homology polytope.   The method of proof 
adopted here, which involves  
only elementary considerations,  is  
quite different in spirit and applies equally well to
 previously established cases.  

\section{$K$-theory of smooth complete toric varieties}
Let $\rho$ be any edge in $\Delta$. 
As $X$ is smooth, the Weil divisor $V(\rho)$ determines a $\bt$-equivariant 
line bundle $\cO(V(\rho))$ which will be denoted $L_\rho.$  The bundle $L_\rho$ 
admits a ($\bt$-equivariant) algebraic cross-section $s_\rho\colon X\lr L_\rho$ 
which vanishes to order $1$ along $V(\rho)$, that is, the zero scheme 
of $s$ equals the variety $V(\rho)$.  (This is a general fact 
concerning the line bundle $\cO(D)$ associated to an effective Weil divisor $D$ on a smooth variety.  Any line bundle can be regarded as a subsheaf of 
$\kappa_X$, the sheaf of  `total quotient rings' of $X$.  When $D$ 
is effective, $\cO(-D)\hookrightarrow \kappa_X$  is 
the ideal sheaf of  the subscheme $D$. 
See \cite[Ch. II, Proposition 6.18]{hart}. The inclusion $\cO(-D)\hookrightarrow \cO_X$ 
yields, upon taking duals, a morphism $s_D\colon \cO_X\lr \cO(D)$ which 
defines the required global section with scheme $D$.  
Applying these considerations to $L_\rho$ we obtain the required 
section $s_\rho.$  Furthermore,  $c_1(L_\rho)\in A^1(X)$ equals the divisor class 
$[V(\rho)]$.  Topologically, the first Chern class $c_1(L_\rho)\in H^2(X)$ 
can also be described as  the cohomology class -- also denoted $[V(\rho)]\in H^2(X)$ -- dual to the submanifold $V(\rho)\subset X$, 
which is the image of $[V(\rho)]\in A^1(X)$ under the isomorphism $A^*(X)\cong H^*(X)$. 

If $\rho_1,\cdots, \rho_k$ 
are edges of $\Delta$ which do not span a cone of $\Delta$, then the 
section $s\colon X\lr L_{\rho_1}\oplus\cdots\oplus L_{\rho_k}$ defined as 
$s(x)=(s_{\rho_1}(x),\cdots, s_{\rho_k}(x))$ is nowhere vanishing since 
$V(\rho_1)\cap \cdots\cap V(\rho_k)=\emptyset$. Therefore $s$ defines 
a monomorphism 
$\wt{s}\colon \cO\lr L_{\rho_1}\oplus\cdots\oplus L_{\rho_k}$ of vector bundles.  
Denote by $E$ the quotient of 
$L_{\rho_1}\oplus\cdots\oplus L_{\rho_k}$ by the image  of $\wt{s}$.  
Applying the $\gamma$-operation in $K(X)$, 
we see that $\gamma^k([L_{\rho_1}\oplus\cdots\oplus L_{\rho_k}]-k)=\lambda^k([E])=0$ as rank of $E$ equals $k-1$.
On the other hand $\gamma^k([L_{\rho_1}\oplus\cdots \oplus L_{\rho_k}]-k)
=\prod_{1\leq i\leq k} \gamma^1([L_{\rho_i}]-1)$.  That is, 
$$\prod_{1\leq i\leq k}([L_{\rho_i}]-1)=0.\eqno{(3)}$$

Let $u\in M$. Set $L_u:=\prod_{\rho\in \Delta(1)}L_\rho^{\langle u,v_\rho\rangle}.$  
The first Chern 
class of the line bundle $L_u$ can be readily calculated to be 
$\sum_{\rho\in \Delta(1)}\langle u,v_\rho\rangle c_1(L_\rho)=\sum 
\langle u,r_\rho\rangle [V(\rho)]=0$ in view of equation (2) above. 
Since the isomorphism class of a line bundle is determined by its 
first Chern class,  we have the following equation in $K(X)$:
$$\prod_{\rho\in \Delta(1)}[L_\rho]^{\langle u,v_\rho\rangle}=1.\eqno{(4)}$$ 

\begin{definition}\label{cand}  Let $\Delta$ be a complete regular fan in $N$.
Let $R(\Delta)$ denote the ring $\bz[x_\rho\mid \rho\in \Delta(1)]/\frak{I}$ where 
$\frak{I}$ is ideal generated by the elements:\\
(i) $x_{\rho_1}\cdots x_{\rho_k}=0$ whenever $\rho_1,\cdots, \rho_k\in \Delta(1)$ 
do not span a cone of $\Delta$. \\
(ii) $z_u:=\prod_{\rho\in \Delta(1)\mid\langle  u,v_\rho\rangle>0}(1-x_\rho)^{\langle u,v_\rho\rangle}-\prod_{\rho\in \Delta(1)\mid\langle u,v_\rho\rangle<0}
(1-x_\rho)^{-\langle u,v_\rho\rangle}$.
\end{definition}

We are now ready to state the main theorem.

\begin{theorem}\label{main}
The $K$-ring of the complete non-singular toric variety $X(\Delta)$ 
is isomorphic to $R(\Delta)$ under an isomorphism $\psi$  which maps 
$x_\rho$ to $(1-[L_\rho])$. 
\end{theorem}
  
In view of equations (3) and (4) it is clear that there is a homomorphism 
of rings $\psi\colon R(\Delta)\lr K(X(\Delta))$. In \S3 we shall 
show that $\psi$ is onto. In \S4 we complete the proof by showing 
that both the abelian groups $R(\Delta)$ and $K(X(\Delta))$ are 
free of the same rank. In \S5 we shall extend the results to 
the class of torus manifolds with orbit space a homology polytope. 

\section{Line bundles and $K$-theory}
In this section $X$ denotes a path connected finite CW complex. 
Assume that $H^*(X)$ is generated by $H^2(X)$.  Then the $K$-ring 
of $X$ is generated (as a ring) by the classes of line bundles on $X$.  
This was proved in \cite{su1} under the hypothesis that $X$ has cells only 
in even dimensions. However essentially the same proof 
works under our weaker hypothesis. Indeed,  
suppose that $\dim(X)<2n$ and that $H^2(X)$ is generated as an abelian 
group by $k$ elements.  Then one has a continuous map 
$f\colon X\lr (\bc\bp^n)^k$  which induces a surjection in cohomology in  
dimension $2$ and hence, by our hypothesis on $X$, in all dimensions. 
Since the cohomology of $X$ vanishes in odd dimensions, it follows 
that the Atiyah-Hirzebruch sequence collapses. Since $f^*$ induces 
surjection in cohomology, the naturality of 
the spectral sequence implies that it $f^*$ induces surjection in 
$K$-theory.  Since $K(\bc\bp^n)$ is generated by line bundles, it 
follows by the K\"unneth theorem for $K$-theory that $K((\bc\bp^n)^k)$ 
is generated by line bundles.  Hence $K(X)$ is also generated by 
line bundles.  As a consequence we obtain 

\begin{proposition} 
Let $X(\Delta)$ be a complete non-singular toric variety. 
With notations as in \S 1, 
the ring homomorphism $\psi\colon R(\Delta) \lr K(X(\Delta))$ is a surjection.
\hfill $\Box$
\end{proposition}

\begin{remark}\label{lowerbound}
Suppose that $H^*(X)=H^\ev(X)$ is a free abelian group. A  
straightforward argument involving the 
Atiyah-Hirzebruch spectral sequence shows that, as far as the additive structure 
is concerned, $K(X)$ is a free abelian group of 
rank equal to the Euler characteristic $\chi(X)$ of $X$.  

It well-known that $\chi(X(\Delta))$ equals $\# \Delta(n)$, the 
number of $n$-dimensional cones in $\Delta$ \cite{fulton}. 
Since $\psi$ is a surjection we conclude that as an abelian 
group, the rank of $R$ is at least $\#\Delta(n)$.
\end{remark}

\section{Proof of  Theorem \ref{main}} 
Let $\Delta$ be a complete regular fan in $N$ and let $R(\Delta)$ be 
the ring defined in \S1. It is clear that the set of all monomials 
$x(\sigma):=x_{\rho_1}\cdots x_{\rho_k}, \sigma\in \Delta,$ where 
$\sigma\in \Delta(k)$ is spanned by edges $\rho_1,\cdots,\rho_k$ 
forms a generating set for $R(\Delta)$. We shall show that $R(\Delta)$ 
is a free abelian group and describe a monomial basis for it.

\noindent{\bf Filtered rings}\\
Let $S$ be the polynomial ring $\bz[x_1,\cdots, x_d]$ and let $\ffi\subset S$ 
be an ideal generated by elements $f_1,\cdots,f_m$ where each $f_j$ has 
constant term zero. The ring $S$ is graded where we set $\deg(x_i)=1$. 
We denote the abelian group of all homogeneous polynomials of degree 
$j$ by $S_{(j)}$.  Consider the multiplicative filtration 
$S=S_0\supset S_1\supset \cdots\supset S_r\supset \cdots$ where 
$S_r=\oplus_{j\geq r}S_{(j)}.$ 
In view of our assumption about the generators of $\ffi$, 
this is an $\ffi$-filtration, i.e., $\ffi S_r\subset S_{r+1}$ for $r\geq 1$.  
Let $R$ denote the quotient ring $S/\ffi$. The filtration of $S$ 
induces a decreasing multiplicative 
filtration $R_0\supset R_1\supset \cdots $ of $R$. 
Let $\gr(R)$ denote the associated graded ring. 
The group of all homogeneous element of degree $j$ in $\gr(R)$ 
is denoted $\gr(R)_{(j)}$.  Clearly 
$\gr(R)_{(j)}=R_j/R_{j+1}=S_j/(\ffi\cap S_{j}+S_{j+1})\cong (S_j/S_{j+1})/
(\ffi\cap S_j/S_{j+1})$ for $j\geq 0$.

Suppose that $\ci$ is the ideal of $S$ generated by the initial forms of 
$f_1,\cdots, f_m$. (By the initial form of $f\in S_r\setminus S_{r+1}$ 
we mean the homogeneous polynomial $\inn(f)$ of degree 
$r\geq 1$ where $f-\inn(f)\in S_{r+1}.$)
Set $\wt{R}:=S/\ci$. Since $\ci$ is a homogeneous ideal, the ring $\wt{R}$ 
inherits a grading from $S$. 
The group $\wt{R}_{(j)}$ of all homogeneous elements of degree $j$ 
is isomorphic to $S_{(j)}/\ci\cap S_{(j)}\cong 
(S_j/S_{j+1})/((\ci\cap S_j+S_{j+1})/S_{j+1})$.
Note that $\ci\cap S_{j}+S_{j+1}\subset \ffi\cap S_j+S_{j+1}.$ 
Hence there exists a natural epimorphism $\eta_j\colon\wt{R}_{(j)}\lr \gr(R)_{(j)}$. 
We let $\eta\colon \gr(\wt{R})\lr \gr(R)$ be direct sum of $\eta_j, j\geq 0$.

Now let $S=\bz[x_\rho\mid\rho\in \Delta]$ and let the ideal $\ffi$ and the set of 
generators of $\ffi$ be as in Definition \ref{cand}. 
The initial form of $z_u$ is seen to be 
$\inn(z_u)=\sum_{\rho\in \Delta(1)}\langle u,v_\rho\rangle x_\rho=:h_u$. 
Thus the ideal $\ci$ is generated by the set of monomials listed in 
\ref{cand} (i) and the elements $h_u,~u\in M.$  By \cite[Theorem 10.8]{danilov}
we see that the ring $\wt{R}(\Delta)=S/\ci$ is isomorphic to the Chow ring  
$A^*(X(\Delta))$ under an isomorphism which maps $x_\rho$ to $[V(\rho)]$. 
We shall identify $\wt{R}(\Delta)$ with $A^*(X(\Delta))$. Thus we 
obtain a surjective homomorphism of graded abelian groups $\eta\colon 
A^*(X(\Delta))\lr \gr(R(\Delta))$.  

\indent\textsc{Proof of Theorem \ref{main}.} We shall abbreviate $X(\Delta)$ to $X$ etc.
Recall that $A^*(X)\cong H^*(X)$ is a free abelian group, its 
rank being equal to $\#\Delta(n)$, the number of 
$n$-dimensional cones of $\Delta$. Since $\eta$ is a surjection, 
$\gr(R(\Delta))$ is an abelian group generated by at most $\#\Delta(n)$  
elements.  Since $R$ and $\gr(R)$ 
have equal rank, it follows that $R$ must be 
free abelian of rank {\it at most} $\#\Delta(n)$. In view 
of Remark \ref{lowerbound}, we conclude that rank of $R$ {\it 
equals} $\#\Delta$.   

Since $\psi\colon R\lr K(X)$ is a surjection between free abelian 
groups of same rank, it follows that $\psi$ is an isomorphism.
This completes the proof. \hfill $\Box$

Let $X$ be a smooth complete variety over $\bc$.
Consider the `forgetful' homomorphism $f\colon\ck(X)\lr K(X)$ where 
$\ck(X)$ denotes the Grothendieck $K$-ring of algebraic vector bundles 
over $X$.  Since $X$ is smooth, $\ck(X)$ is isomorphic to the Grothendieck 
group $\ck'(X)$ of coherent sheaves over $X$. 
One has a `topological filtration' on $\ck'(X)\cong\ck(X)$ and one has 
a well-defined surjective homomorphism $\phi\colon A_*(X) \lr \gr(\ck(X))$
of graded abelian groups.  When $X=X(\Delta)$, it follows that 
$\gr(\ck(X))$  is a quotient of the free abelian group $A_*(X)=A^{n-*}(X)$. 
It follows that $\ck(X)$ is generated by at most $\chi(X)$ elements.
Since $K(X)$ has rank equal to $\chi(X)$, it follows that $f\colon 
\ck(X)$ is an isomorphism of rings. We record this as 

\begin{theorem}
Let $X=X(\Delta)$ be a smooth complete toric variety. Then, 
the forgetful morphism of rings $f\colon \ck(X)\lr K(X)$ is an isomorphism.
In particular, $\ck(X(\Delta))$ is isomorphic to $R(\Delta)$.  \hfill $\Box$
\end{theorem}

\begin{remark}
(i) It is immediate from the main theorem that $x_\rho\in R$ are nilpotent since 
$[L_\rho]$ are invertible.  
Indeed, the surjectivity of $\eta$ already implies that the filtration 
$R_j=R_{n+1}$ for all $j>n$. Since $\cap_{j\geq 0} R_j=0$, it follows that 
$R_j=0~\forall j>n$.  \\
(ii) It follows from our proof that $\eta\colon A^*(X)\lr\gr(R(\Delta))$ is 
an isomorphism.  Recall that if $\sigma\in \Delta(k)$ is spanned by edges 
$\rho_1,\cdots, \rho_k$, then $[V(\sigma)]=[V(\rho_1)]\cdots [V(\rho_k)]\in A^{n-k}(X)$.  If $\Gamma\subset \Delta$ is a collection of cones 
such that $\{[V(\sigma)]\}_{\sigma\in \Gamma} $ is $\bz$-basis for $A^*(X)$, then 
$\{x(\sigma) \}_{\sigma\in \Gamma}$ is a monomial basis for $R(\Delta)$. \\
\end{remark}

\section{K-theory of torus manifolds}
The algebraic geometric notion a non-singular projective toric variety
has been brought home to realm of topology by  
Davis and Januszkiewicz  \cite{dj} who introduced the class of 
quasi-toric manifolds. (Davis and Januszkiewicz called them `toric manifolds'; 
the present terminology of  `quasi-toric manifolds' is due to Buchstaber and 
Panov \cite{bp}.) See also Masuda \cite{masuda}.  
Recently, Masuda and Panov \cite{mp} introduced a new class of manifolds called 
{\it torus manifolds} in their efforts to develop the analogue of a non-singular 
complete toric variety. The class of torus manifolds, which includes 
all complete non-singular complex toric varieties, is much more 
general than that of quasi-toric manfold as  
there are non-singular complete toric varieties which are not quasi-toric 
manifolds.  

In this section we obtain a description of the $K$-ring of  a torus manifold $X$ 
assuming that the orbit space $X/T=:Q$ is a {\it homology polytope}.  
Most of the ingredients needed for the proof of Theorem \ref{torus} 
can be found in \cite{mp}.  

We shall recall here the results of \cite{mp} 
concerning the cohomology of toric manifolds with locally standard action and 
orbit space a homology polytope that are needed for our purposes. 
We recall below the definition and some basic facts relevant for our purposes,  
referring the reader to \cite{mp} for a detailed exposition of torus manifolds. 
As mentioned in the Introduction, our main result in this section, namely, Theorem \ref{torus} 
subsumes Theorem \ref{main}.  In fact, as we shall see, 
the steps involved in the proof in both cases 
are the same; we need only to establish those 
steps whose proofs in the algebraic geometric situation  
are either unavailable or are not obvious in the topological context of 
torus manifolds. 

We begin by recalling the basic definition and properties of torus manifolds. 
Let $T:=(\bs^1)^n$ denote the {\it compact} torus and let $X$ be a smooth 
compact oriented connected manifold of dimension $2n$ on which $T$ acts 
effectively with a finite non-empty set of $T$-fixed points. Such a manifold $X$ 
is called a torus manifold. The orbit space $Q=:X/T$ is a `manifold with 
corners'. The $T$-action on $X$ is called {\it locally standard} if  $X$ is covered 
by $T$-invariant open sets $U$ such that $U$ is equivariantly diffeomorphic 
to an invariant open subset contained in a $T$-representation 
$\mathcal{U}\cong \bc^n$ 
whose characters form a $\bz$-basis for $\Hom(T,\bs^1)\cong \bz^n$. 
For example, $X=X(\Delta)$ a complete non-singular toric variety with action 
of $\bt$ restricted to $T\subset \bt$ is a torus manifold with locally standard 
$T$-action.  Also, any quasi-toric manifold  
is a torus manifold. A characteristic submanifold of $X$ is a codimension $2$ submanifold which is pointwise fixed by a one-dimensional subgroup of $T$.  
There are only finitely many characteristic submanifolds; we denote them by
$V_1,\cdots, V_d$.  

In the case when $X$ is a smooth projective toric variety or a quasi-toric manifold 
the orbit space $Q$ is a simple polytope.  However, for a general 
torus manifold, this is not so. Assume that $X$ is torus manifold with locally 
standard $T$ action. 
Define the boundary of $Q$, denoted $\partial Q$, to be the set of all points which do {\it not} have a neighourhood homeomorphic to an open set of $\br^n$. Then $\partial Q$ is the image under the quotient map $\pi\colon X\lr Q$ of the union of all characteristic 
submanifolds of $X$.  Denote the image of $V_i$ by $Q_i\subset \partial Q, 1\leq i
\leq d$; these are the {\it facets} of $Q$. A non-empty intersection of 
facets are called {\it prefaces} of $Q$; a {\it face} of $Q$ is a connected 
component of a preface. The space $Q$ is called {\it homology 
polytope} if every preface $F$, including $Q$ itself, is acyclic, i.e., $\wt{H}_*(F;\bz)=0$,  in particular, every preface of $Q$ is path connected. One of the 
main results of \cite{mp} is that the integral cohomology of $X$ is generated 
by degree $2$ elements if and only if the action is locally standard and $Q$ 
is a homology polytope.   
{\it In this paper we shall only consider torus manifolds 
with locally standard $T$-action whose orbit space is a homology polytope.}

Note that any characteristic submanifold of such a torus manifold inherits 
these properties; see \cite[Lemma 2.3]{mp}.  In particular, characteristic 
submanifolds are orientable and have $T$-fixed points. We fix an omni-orientation 
of $X$, i.e., orientations of $X$ and all its characteristic submanifolds. Thus 
the normal bundle $\nu_i$ to the imbedding $V_i\hookrightarrow X$ is 
oriented by the requirement that $\nu_i\oplus TV_i$ be oriented-isomorphic 
to $TX|V$, the tangent bundle of $X$ restricted to $V$.

\noindent{\bf Line bundles over $X$}\\
Let $V\subset X$ a codimension $2$ oriented closed connected submanifold of an 
{\it arbitrary} oriented closed connected manifold of dimension $m$.  
Let $[V]\in H^2(X;\bz)$ denote the cohomology class dual to $V$ and let $L$ denote a complex line bundle over $X$ with $c_1(L)=[V]$.   Assume that $H^1(V;\bz)=0$. 

\begin{lemma}\label{section} Suppose that $c_1(L)=[V]\neq 0$.  Then the complex line bundle $L$ admits a section $s\colon X\lr L$ such that the zero locus of $s$ is precisely $V$.
\end{lemma}
\noindent
\begin{proof} Assume that $X$ has been endowed with a Riemannian metric and 
let $N\subset X$ denote  a tubular neighbourhood of $V$ 
which is identified with the disk bundle associated to  
the normal bundle $\nu$  to the imbedding $V\hookrightarrow X$. The 
normal bundle is canonically oriented by the requirement that $TV\oplus \nu$ 
is orientated isomorphic to $TX|V$.  
Since $\nu$ is an oriented real $2$-plane bundle, we regard it as a complex line 
bundle. We shall denote by $\pi\colon N\lr V$ the projection of the disk bundle. 

The complex line bundle $\pi^*(\nu)$ over $N$ 
evidently admits a canonical cross-section $\eta\colon N\lr \pi^*(\nu)$ which vanishes 
precisely along $V$.   Take the trivial complex line bundle $\varepsilon$ over  
$(X\setminus int(N))$ and consider the vector bundle map $\wt{\eta}\colon 
 \varepsilon|\partial N\lr \pi^*(\nu)|\partial N$ defined by  $\wt{\eta}(x,1)
 =\eta(x)$ for all $x\in \partial N $. 
Gluing  $\varepsilon|\partial N$ along $\pi^*(\nu)|\partial N$ using the 
identification $\wt{\eta}$ yields a vector bundle $\xi$ over $X$ which 
clearly admits a cross-section $s$ 
(which restricts to $\eta$ on $N$ and $x\mapsto (x,1)$ on $X\setminus  int(N)$) 
that vanishes precisely along $V$.  Now it suffices to 
to establish that $\xi$ is isomorphic to $L$. 
 
Both $\xi$ and $L$ restrict to trivial bundles over $X\setminus int(N)$.
The bundle $L$ restricts to the normal bundle $\nu$ over $V$ (cf. 
\cite[Theorem 11.3]{ms}).  Since $\pi\colon N\lr V$ and 
$V\hookrightarrow N$ are homotopy inverses, it follows that $\xi$ and $L$ 
restrict to isomorphic bundles over $N$ as well.  Let $\sigma\colon  
\varepsilon|\partial N\lr L|\partial N\cong \pi^*(\nu)|\partial N$ 
be the gluing data for $L$.   The vector bundle maps $\sigma $ and $\wt{\eta}$ 
are related by a map $\partial N\lr GL_1(\bc)$.  The  homotopy 
classes of such 
maps is isomorphic to $H^1(\partial N;\bz)$. Using the Serre spectral sequence  
of the circle bundle $\partial N\lr V$ hypotheses 
that $H^1(V)=0$ and $c_1(L)\neq 0$ we see that $E_3^{0,1}=0, E_2^{1,0}=0$. 
Therefore $H^1(\partial N;\bz)=0$.  Hence $\sigma$ and 
$\wt{\eta}$ are homotopic via bundle isomorphisms. It follows that $L$ is isomorphic  to $\xi$. 
\end{proof}

Now applying the above lemma to characteristic submanifolds $V_1,\cdots, 
V_d$ of a torus manifold $X$  with locally standard $T$-action and orbit 
space a homology polytope, we see that there exist complex line bundles $L_1,\cdots, 
L_d$ such that $c_1(L_i)=[V_i]$ and each $L_i$ admits a section $s_i\colon X
\lr L_i$ which vanishes precisely along $V_i$, $1\leq i\leq d$. 
We proceed as we did in obtaining Equation (3) in \S2, to 
obtain the following equation in $K(X)$  
$$\prod_{1\leq j\leq r}(1-[L_{i_k}])=0 \eqno{(5)}$$ 
whenever $\cap_{1\leq k\leq r}V_{i_k}=\emptyset$.

\noindent{\bf The characteristic map}\\
Recall that 
the characteristic submanifolds $V_i$ are, by definition, fixed by one dimensional 
subgroups $S_i$ of $T$. Our assumption on $X$ (local standardness and 
orbit space being a homology polytope) implies that every characteristic 
submanifold has a $T$-fixed point. There is unique $1$-parameter subgroup $v_i\in \Hom(\bs^1,T)\cong \bz^n$ by the following requirements, the first of which determines 
$v_i$ upto sign: (i) $v_i$ is a primitive element in $\Hom(\bs^1,T)$ with 
image $S_i$, and,  
(ii) the sign ($+v_i$ or $-v_i$) is determined by orienting $S_i$ so 
that at any point $p\in V$,  the oriented normal plane $\nu_p$ 
is {\it oriented} isomorphic to the tangent space to $S_i$ at the identity element.  

The map $\Lambda\colon\{Q_1,\cdots,Q_d\}\lr \Hom(\bs^1,T)$ which 
maps $Q_i$ to $v_i$ is called the {\it characteristic map}.  Under our hypothesis 
of local standardness and $Q$ being a homology polytope, 
the manifold $X$ is determined upto equivariant 
diffeomorphism by the pair $(Q,\Lambda)$.  (See \cite[Lemma 4.5]{mp}.) 

Local standardness of the $T$-action implies that if $V_{i_1}\cap\cdots\cap V_{i_r}
\neq \emptyset$, then $\{v_{i_1},\cdots,v_{i_r}\}$ is part of a $\bz$ basis for 
$\Hom(\bs^1;T)\cong \bz^n$.

\noindent{\bf Cohomology of $X$} \\
We now recall from \cite[Corollary 7.8]{mp} the description of the integral cohomology ring of $X$. 

Let $\mathcal{Q} $ denote the set of all faces of $Q$. Then 
the cohomology of $X$ has the following description:\\
$H^*(X;\bz)\cong \bz[x_F;F\in \mathcal{Q}]/I$ where $I$ is the ideal generated by 
the following two types of elements:\\
(i) $x_Ax_B-x_{A\vee B}x_{A\cap B}$ where $A\vee B$ denotes the 
smallest face of $Q$ which contains both $A$ and $B$,\\
(ii) $\sum_{1\leq i\leq  d}\langle u,v_i\rangle x_{Q_i},~ u\in \Hom(T,\bs^1)$ 
where $v_i=\Lambda(Q_i)\in \Hom(\bs^1;T).$\\  The element $x_{Q_i}$ 
corresponds under the isomorphism to $[V_i]\in H^2(X;\bz).$
(It is understood that $x_Q=1,$ and $x_\emptyset=0$.)

From our hypothesis, $H^*(X;\bz)$ is generated by degree two elements. 
Set $x_i:=x_{Q_i}$, for $1\leq i\leq d$.
Any face $F$ of $Q$ is the intersection of those facets of $Q$ which contain $F$.
If $F$ is of codimension $r$, then it is contained in exactly $r$ 
distinct facets, say, $Q_{i_1},\cdots , Q_{i_r}$ and so the intersection 
$F=Q_{i_1}\cap\cdots\cap Q_{i_r}$ is transversal. Therefore 
$x_F=x_{i_1},\cdots,x_{i_r}$. Thus, we see that 
$H^*(X;\bz)=\bz[x_i;1\leq d]/\ci$ where the ideal $\ci$ of relations is 
generated by the elements:\\
(iii) $x_{i_1}\cdots x_{i_r}$ whenever $V_{i_1}\cap\cdots\cap V_{i_r}=\emptyset$\\
(iv) $\sum_{1\leq i\leq d}\langle u,v_i\rangle  x_i$.

\begin{proposition}\label{toruscohomology}(Masuda-Panov \cite{mp})\\
Let $X$ be a $T$-torus manifold with locally standard action whose orbit space is a homology polytope. With the above notations, 
$H^*(X;\bz)\cong \bz[x_1,\cdots,x_d]/\ci$ where $x_i$ represents 
the cohomology class dual to the characteristic manifolds 
$V_i, 1\leq i\leq d$.  In particular, $H^*(X;\bz)$ is a free abelian group. 
\end{proposition}
\noindent
\begin{proof}The only part of the theorem that needs to be established is 
that $H^*(X;\bz)$ is a free abelian group as other assertions follow 
from \cite[Corollary 7.8]{mp} as noted above.
By \cite[Theorem 7.7]{mp}, the cohomology of $X$ with 
coefficients in $\bz$ or $\bz/p\bz$ for any prime $p$ vanishes in 
odd dimensions, it follows that $H^*(X;\bz)$ is torsion free. As $X$ is 
a compact manifold, $H^*(X;\bz)$ is finitely generated. It follows that 
$H^*(X;\bz)$ is free abelian.
\end{proof}

Recall that $L_i$ is the complex line bundle over $X$ with $c_1(L_i)=[V_i]$.
As an immediate consequence of the above description of $H^*(X;\bz)$, 
we obtain the following equality in $K(X)$, analogous to equation (4):\\
$$\prod_{1\leq i\leq d}[L_i]^{\langle u,v_i\rangle}=1.\eqno(6)$$

We are now ready to state the main result of this section.  

\begin{theorem}\label{torus} 
Let $X$ be any $T$-torus manifold with locally standard action whose 
orbit space is $Q$ a homology polytope. Then $K(X)$ is the ring isomorphic to 
the ring 
$R(Q,\Lambda):=\bz[y_1,\cdots,y_d]/\fj$ where $\fj$ is the ideal generated by the 
elements\\
(i) $y_{i_1}\cdots y_{i_r}$ whenever $Q_{i_1}\cap\cdots\cap Q_{i_r}=\emptyset$\\
(ii) $\prod _{\{i\leq d: \langle u,v_i\rangle>0\}}(1-y_i)^{\langle u,v_i\rangle}-\prod_{
\{j\leq d:\langle u,v_j\rangle<0\}}(1-y_j)^{-\langle u,v_j\rangle}$ for each $u\in \Hom(T;\bs^1),$ where $v_i:=\Lambda(Q_i), 1\leq i\leq d$. The isomorphism   
is established by sending $y_i$ to $1-[L_i]$. 
\end{theorem}
\begin{proof}
The proof follows exactly as in the case of non-
singular complex toric varieties. 
In view of Proposition \ref{toruscohomology},  $K(X)$ is generated 
by $[L_i], 1\leq i\leq d$; see \S3. Furthermore, $K(X)$ is a free abelian 
group of rank $\chi(X)$.  Set $R:=R(Q,\Lambda)$. 
From equations (5) and (6), there is a ring homomorphism $\psi\colon 
R\lr K(X)$ which maps $y_i$ to $1-[L_i]$ which is {\it surjective}.
Arguing as in \S4, we see that there is a decreasing filtration on $R$ 
and a surjective homomorphism of abelian groups $H^*(X;\bz)\lr 
\gr(R)$, showing that  as an abelian group, $R$  
generated by at most $\chi(X)$ many elements. 
Since $\psi$ is surjective, it follows that $R$ also a free abelian 
group of rank $\chi(X)$ and  that $\psi$ is an isomorphism of {\it rings}. 
\end{proof}

\end{document}